%% file: BMCircle2.tex
\documentclass[12pt, english]{article}
\input{figure_preamble.tex}
\usepackage{amssymb,amsfonts,amsmath,amsthm}
\usepackage{lineno,hyperref}
\usepackage{bbm}
\usepackage{dsfont}
\usepackage{graphicx}
\usepackage[margin=1in]{geometry}
\usepackage{enumitem}
\modulolinenumbers[5]

\newcommand{\beqn}{\vspace{-0.25cm}\begin{eqnarray*}}
\newcommand{\eeqn}{\end{eqnarray*}}
\newcommand{\bneqn}{\vspace{-0.25cm}\begin{eqnarray}}
\newcommand{\eneqn}{\end{eqnarray}}

\renewcommand{\exp}[1]{\mathrm{exp}\parens{#1}}

\newcommand{\indic}[1]{\mathds{1}_{#1}}
\newcommand{\beq}{\begin{eqnarray*}}
\newcommand{\feq}{\end{eqnarray*}}
\newcommand{\feqn}{\end{eqnarray}}

\newcommand{\bracks}[1]{\left[#1\right]}
\newcommand{\parens}[1]{\left(#1\right)}
\newcommand{\expe}[1]{\mathbb{E}\bracks{#1}}
\newcommand{\natlog}[1]{\ln\parens{#1}}
\newcommand{\prob}[1]{\mathbb{P}\parens{#1}}
\newcommand{\mathand}{~~\text{and}~~}

\newcommand{\inverse}[1]{\parens{#1}^{-1}}
\newtheorem{theorem}{Theorem}
\makeatletter \@addtoreset{theorem}{section}\makeatother

\newtheorem*{theorem*}{Theorem}

\newtheorem{remark}[theorem]{Remark}

\title{On the Time for Brownian Motion to Visit Every Point on a Circle}
\author{Philip Ernst\footnote{Rice University} \, and Larry Shepp \footnote{The Wharton School of the University of Pennsylvania (Deceased April 23, 2013)}}
\date{}
\begin{document}
\maketitle

\begin{abstract}
Consider a Wiener process $W$ on a circle of circumference $L$. We prove the rather surprising result that the Laplace transform of the distribution of the first time, $\theta_L$, when the Wiener process has visited every point of the circle can be solved in closed form using a continuous recurrence approach. 
\end{abstract}


\section{Introduction}
Consider a Wiener process on a circle of circumference
$L$. The distribution of the first time, $\theta_L$, when the Wiener process
has visited every point of the circle is equivalent, via the natural
bijection between and interval of the form $[b, b+L)$ on the real line and a circle of circumference $L$, to the distribution of the first time when the range of the Wiener process on the real line is of length $L$. This distribution is well-known and it has the following Laplace transform: [see, for example, ~\cite{Borodin}, p.242]

\bneqn \label{central}
\expe{e^{-s\theta_L}} = \frac{1}{\cosh^2 \parens{L\sqrt{\frac{s}{2}}}}, ~~s \ge 0.
\eneqn
~\cite{Feller}, in writing about the range of a Wiener process, did so using explicit probability density calculations. ~\cite{Imhof} discovered Laplace transform for the first time, $\theta_L$, when the Wiener process has visited every point of the circle,  again via explicit probability density calculations.  Further computations employing the Laplace transform for $\theta_L$ were presented in ~\cite{Vallois}. However, in departure from these previous works, we prove the result in equation (\ref{central}) using a continuous recurrence setup. We do so by calculating the left hand side in terms of random
variables representing how long it takes to cover a range of length
$L$, \emph{given that one is already at an endpoint of a range of length $a$}
(which counts as being covered already).  This is the idea behind the
definition of $\theta_{a, L}$, which is defined in Section \ref{recurrence}.

Key to our recurrence will be the concept of a \emph{switchback}.  
Imagine we pick some $a \in \mathbb{R}^{+}$ that is less than $L$.
Consider the maximum, $M_a$, of $W$ until the first visit to the point, $-a$ on the negative half-axis. (Here, $M_a>0$; otherwise, the process would have moved directly from 0 to $-a$, which occurs w.p. 0) We call the time of this first visit $\tau_{-a}$. We say that a 
``{\em switchback}'' occurs when $W$ hits $-a$ before the length of the
range, $a+M_a$, is $L$. Formally, let $\indic{a,L}$ be the indicator random variable for the event of a switchback, defined as follows:

\beqn
\indic{a,L}=
\begin{cases}
1 \quad \text{if} \,\,\, \text{inf}  \{ t: 0 \leq t < \infty\, | \,W_t=-a \}\leq \text{inf}  \{ t: 0 \leq t < \infty\, | \,W_t=L-a \}\\ 
0 \quad \textit{otherwise}.
\end{cases}
\eeqn

   After a switchback, the process continues from $-a$
with a starting range of $M_a+a$ (i.e., the interval $[-a, M_a]$ has
been covered). By translation and reflection invariance, as well as the symmetry of Brownian motion, we may just as
well assume that we are at the point $0$ and have covered the
interval $[-(a+M_a), 0]$.  We then repeat the process and say that a
second switchback occurs if we reach $-(a+M_a)$ before covering a
range of length $L$.  To summarize:\\

\noindent \textit{Step 1}: We start our process at the right hand end
of $[-a,0]$ and we consider this interval as already being covered. $M_a$ is the maximal value attained before the time
$\tau_a$ that we first hit $-a$.  The total range is $a+M_a$. If $M_a
\geq L-a$, then we have covered an interval of length $L$ before
reaching $-a$, and no switchback occurs.  If not, a switchback
occurs and we continue to Step 2.\\\\
\textit{Step 2}: We have covered a range of length $a+M_a$. Without loss of generality, we consider the interval $[-(a+M_a),0]$ to have been covered. Let $-(a+M_a)$:=$-a^\prime$, and
start the process on the right hand end of $[-a^\prime, 0]$. If
$M_{a^\prime} \geq L-M_{a^\prime}$, no switchback occurs. Otherwise, another
switchback occurs and we continue to Step 3.\\\\
\textit{Step 3}: We have covered a range of length $a^\prime+M_{a^\prime}$. Without loss of generality, we consider the interval $[-(a^\prime+M_{a^\prime}),0]$ to have been covered. Let $-(a^\prime+M_{a^\prime})$ be called $-a^{\prime\prime}$, and start the process on the right hand end of $[-a^{\prime\prime}, 0]$. If $M_{a^{\prime\prime}} \geq L-a^{\prime \prime}$, a switchback occurs. Otherwise, continue Step 3 recursively until a range of length $L$ has been covered. \\\\
\noindent Steps 1-3 are illustrated in Figure 1.\\

\begin{figure}[h]
\centering
\input{just_fig.tex}
\caption{Illustration of steps 1-3 in Section 1.}
\label{fig:meyer}
\end{figure}

\noindent In Section 3 we prove that the recurrence can be solved in closed form. In Section 4 we prove that the number $\nu = \nu_{a,L}$, of switchbacks before covering an interval of length $L$ has a Poisson
distribution with parameter $\lambda = \log{\frac{L}{a}}$. Thus, as $a \downarrow 0$, the number of switchbacks goes to infinity at a logarithmic rate. 

\section{Solving the Recurrence}\label{recurrence}

We proceed to solve for the recurrence. First,  consider a Wiener process $W(t), ~t \ge 0$.
For each fixed, $a > 0$, let $M_a$ denote the maximum positive value of $W(t)$
before the first hitting time of $-a$. Assuming that $L-a$ is positive, we have
\beqn
\prob{M_a \le y}= \prob{\tau_{-a} < \tau_{y}} = \frac{y}{a+y},
\eeqn
by the logic of the gambler's ruin.\\

Let $I(t)$ be the range of the Wiener process up to time $t$. Define $\theta_{a,L}$ to be the random variable representing the time until ${I}(t) \cup [-a,0]$ has length $L$. We proceed by defining

\bneqn
f(s,a,L) := \expe{\exp{-s \theta_{a,L}}},
\eneqn

\noindent
where $f(s,a,L)$ is considered a function of $a$ with $s$ and $L$ being held constant. By abuse of notation, we label  $f(s,a,L)$ as $f(a)$.

Let us define the following functions

\bneqn
F(s,y) = \expe{\exp{-s \tau_{-a}} \indic{\tau_{-a} < \tau_y}} \mathand G(s,y) = \expe{\exp{-s \tau_y} \indic{\tau_y < \tau_{-a}}}.
\eneqn

\vskip .2in

\noindent
We now employ the well-known fact (see ~\cite{Borodin}, amongst other sources) that for any $c$,

\begin{equation}\label{martingale}
\exp{cW(t) - \frac{c^2}{2} t}~~ t \ge 0
\end{equation}
is a martingale. If $s = \frac{c^2}{2}$,
we easily obtain the following standard and well known forms of $F(s,y)$ and $G(s,y)$ (see ~\cite{Borodin}, amongst other sources)

\bneqn \label{imp2}
F(s,y) = \frac{\sinh {cy}}{\sinh\parens{{c(a+y)}}}
\mathand
G(s,y) = \frac{\sinh{ca}}{\sinh\parens{{c(a+y)}}}.
\eneqn
Continuing from above, our goal is to write a recurrence for
$f(a)$ in terms of $f(a+y)$ for $0 < y \leq L-a$. To do so, we define $f(a)$ using indicator functions. With the process starting at 0, let the first indicator function represents the case of a switchback, in which $-a$ is hit before the length of the range is $L$. Let the second indicator function denote the case of no switchback.  We may then write

\bneqn
f(a)= \underbrace{\expe{\exp{-s \theta_{a,L}} \indic{\tau_{-a} < \tau_{L-a}}}}_\text{switchback}+ \underbrace{\expe{\exp{-s \theta_{a,L}} \indic{\tau_{L-a} < \tau_{-a}}}}_\text{no switchback}.
\eneqn
Letting $y=L-a$, and using the expression for
$G(s,y)$ in equation $(\ref{imp2})$, we have:

\bneqn
G(s, L-a)= \frac{\sinh{ca}}{\sinh{cL}},
\eneqn
which is exactly the ``no switchback'' term. To calculate the switchback term, we integrate over all possible
values of $M_a$, from $0$ to $L-a$, using $y$ as a dummy variable. $f(a)$ becomes:

\bneqn \label{imp}
f(a) = \frac{\sinh{ca}}{\sinh{cL}} + \int_0^{L-a} f(a+y) \frac{d}{dy} F(s,y) \,dy\,, \quad 0 < a \leq L.
\eneqn

We note that it is possible for $a$ to be $L$ since the original definition of $f(s,a,L)$ gives $f(L)=1$.\\

This recurrence structure is an integral equation. The key idea is that equation
(\ref{imp}) shows that the expected time it takes to get from point
$a$ to point $b$, conditional on starting at a left most point of an
interval, can be found by integrating over all possible left most
points of the subsequent path.

\begin{remark}
This continuous recurrence approach presents an enormously valuable alternative to a direct density calculation. This approach should be helpful in many applied statistical settings in which such a calculation is intractable!
\end{remark}

\section{A Closed Form Solution for the Recurrence}

\begin{theorem}
The recurrence structure in equation (\ref{imp}) can be solved in closed form. Letting $a \downarrow 0$, we obtain $f(0)=\frac{1}{\cosh^2 \parens{L\sqrt{\frac{s}{2}}}}$.
\end{theorem}

\begin{proof}
Using the expression for $F(s,y)$ in equation $(\ref{imp2})$ we have

\bneqn
f(a) = \frac{\sinh{ca}}{\sinh{cL}} + \int_0^{L-a} f(a+y) {\frac{d}{dy} \bracks{\frac{\sinh{cy}}{\sinh{c(a+y)}}}}dy, ~ 0 < a \leq L.
\eneqn

\vskip .2in

\noindent
Differentiating, we obtain:

\beqn
f(a) = \frac{\sinh{ca}}{\sinh{cL}} + \int_0^{L-a} \frac{c\sinh{ca}}{\sinh^2{c(a+y)}} f(a+y)dy, ~0 < a \leq L.
\eeqn

\vskip .2in

\noindent
Finally, substituting $x = a + y$ gives the integral equation

\beqn
f(a) = \frac{\sinh{ca}}{\sinh{cL}} + \int_a^L \frac{c\sinh{ca}}{\sinh^2{cx}} f(x)dx, ~0 < a \leq L.
\eeqn

\vskip .2in

\noindent
Fortunately, this is easy to solve: we divide by $\sinh{ca}$ and let

\beqn
g(x) = \frac{f(x)}{\sinh{cx}}
\eeqn
to arrive at:

\bneqn
g(a) = \frac{1}{\sinh{cL}} + \int_a^L \frac{c}{\sinh{cx}} g(x)dx,~ 0 < a < L.
\eneqn

\noindent
Differentiating with respect to $a$, we obtain the following differential equation for $g$

\bneqn \label{imp3}
g^\prime(a) = -\frac{c}{\sinh{ca}} g(a).
\eneqn

\noindent
Noting that $f(L) = 1$, 

\bneqn \label{imp4}
g(L) = \frac{1}{\sinh{cL}}.
\eneqn

\noindent We now have that

\beqn
g(a) = \frac{1}{\sinh{cL}}\exp{\int_a^L \frac{c}{\sinh{cu}}du},
\eeqn
which is a unique solution to equation $(\ref{imp3})$ with $(\ref{imp4})$ as its initial condition. 

\vskip .2 in

\noindent
Further, 

\beqn
f(a) = \frac{\sinh{ca}}{\sinh{cL}} \exp{\int_a^L \frac{c du}{\sinh{cu}} du}.
\eeqn
We now let $a \downarrow 0$. The limit is

\bneqn
f(0) = \expe{ \exp{-s \theta_L}} = \lim_{a \rightarrow 0} \exp{\int_a^L \parens{\frac{c}{\sinh{cu}} - \frac{c \cosh{cu}}{\sinh{cu}}}du }.
\eneqn

\noindent
Combining the fractions, integrating, and letting $c=\sqrt{2s}$, we obtain

\bneqn
f(0) = \expe{ \exp{-s\theta_L}} = \frac{2}{1+\cosh{cL}} = \frac{1}{\cosh^2 \parens{L\sqrt{\frac{s}{2}}}}.
\eneqn
Since by ~\cite{Oberhettinger}

\bneqn \label{laplace1}
\int_0^\infty \exp{-st} \exp{-\frac{a^2}{2t} }\frac{a}{\sqrt{2\pi}} t^{-\frac{3}{2}}dt = \exp{-a \sqrt{2s}},
\eneqn

\noindent
we can expand $\inverse{\cosh^2{\sqrt{\frac{s}{2}}}}$ in powers of
$e^{-\sqrt{2s}}$ to obtain an infinite series representation of the density
of $\theta_1$. Namely, for $L=1$, we can write:

\bneqn \label{laplace2}
\int_0^\infty e^{-st} p_{_{\theta_1}}(t) dt  \nonumber
&=& 4 ~\exp{-\sqrt{2s}} \parens{1+\exp{-\sqrt{2s}}}^{-2}\\ 
&=& \sum_{n=0}^\infty 4 (-1)^n (n+1) \exp{-(n+1)\sqrt{2s}}.
\eneqn

\noindent
Since the Laplace transform is invertible, it suffices to find a
formula for $p_{_{\theta_1}}(t)$ that makes the above equation true.
If we take

\beqn
p_{_{\theta_1}}(t) = \sum_{n=0}^\infty 4 \frac{(-1)^n (n+1)^2}{\sqrt{2\pi} t^{\frac{3}{2}} } \exp{-\frac{(n+1)^2}{2t}}
\eeqn
and plug this into the left hand side of (\ref{laplace2}), then using
(\ref{laplace1}) gives us the equality (\ref{laplace2}). Since the amount of time to cover a range of length $L$ is equal to $L^2$ multiplied by the amount of time to cover a range of length 1, we write: $\theta_L \sim L^2 \theta_1$. Thus, the density of $\theta_L$ is

\bneqn
p_{_{\theta_L}}(t) = \frac{p_{_{\theta_1}}(\frac{t}{L^2})}{L^2}.
\eneqn
\end{proof}

\section{The Number of Switchbacks}

The result that the number of switchbacks is distributed as a Poisson random variable comes naturally when one accounts for the Markov and scaling properties of Brownian motion. The formal proof of the result follows.

\begin{theorem}
The number of switchbacks
has a Poisson distribution with parameter $\lambda = \log\parens{\frac{L}{a}}.$
\end{theorem}

\begin{proof}
The argument that the number of switchbacks
has a Poisson distribution with parameter $\lambda = \log\parens{\frac{L}{a}}$ is
similar to the argument in Section 3 used to obtain the distribution of $\theta_L$. We begin by defining 

\bneqn
f(a) = f(a,L,t) = \expe{ t^{\nu_{L,a}}},
\eneqn

\noindent where $t$ is a dummy variable and $\nu_{L,a}$ is the number
of switchbacks starting from an endpoint of an interval of length $a > 0$
before the interval grows to length $L$. $\nu_{L,a} = 0$ with probability
$\frac{a}{L}$ and 

\beqn 
\prob{M_a \le y}= \frac{y}{a+y}.
\eeqn
We rewrite $f(a)$ as we did in Section 3, splitting it up by indicator random variables which account for whether or not a switchback has occurred. If no switchbacks have occurred, $\nu_{L,a}=0$, and thus:

\bneqn\label{eqn3}
\expe{ t^{\nu_{L,a}}\indic{M_a>L-a}}=\expe{\indic{M_a>L-a}}=\frac{a}{L}.
\eneqn
Employing equation (\ref{eqn3}) and recalling that a switchback occurs when $M_a< L-a$, we write:

\beqn
f(a) =  \frac{a}{L} + \int_0^{L-a}t^{\nu_{L,y+a}} \prob{M_a \in dy} dy,
\eeqn
which simplifies to

\beqn
f(a) = \frac{a}{L} + \int_0^{L-a}t f(a+y) \frac{d}{dy} \bracks{ \frac{y}{a+y}} dy.
\eeqn

\noindent
Thus, 

\bneqn
f(a) = \frac{a}{L} + t \int_0^{L-a} f(a+y) \frac{a}{(a+y)^2}dy = \frac{a}{L} + at \int_0^L \frac{f(x)}{x^2}dx.
\eneqn

\vskip .2in

\noindent
It is straightforward to verify directly, or to use the differential equation argument used
above, to find the Laplace transform of $\theta_L$, so that

\beqn
f(a) = \exp{\natlog{\frac{L}{a}}(t-1)}.
\eeqn

\noindent Since the above is the form of the Laplace transform of the Poisson distribution, we have the desired result that the number of switchbacks is distributed Poisson with mean $\lambda =\natlog{\frac{L}{a}}$.

\end{proof}

\noindent \textbf{Acknowledgments}: We are greatly appreciative of the two anonymous referees and the Associate Editor whose remarks greatly improved the quality of this manuscript.

\bibliography{mybibfile}

\end{document}

%% file: figure_preamble.tex
\usepackage{tikz}
\usetikzlibrary{snakes}

%% file: just_fig.tex
\tikzset{
  prob_squiggle/.style={decorate, decoration={snake}, draw=black}
}

\newcommand{\lineann}[4][0.5]{%
    \begin{scope}[rotate=#2, black,inner sep=2pt]
        \draw[dashed, black!40] (0,0cm) -- +(0,#1cm)
            node [coordinate] (a) {};
        \draw[dashed, black!40] (#3,0cm) -- +(0,#1cm)
            node [coordinate] (b) {};
        \draw[|<->|] (a) -- node[fill=white] {#4} (b);
    \end{scope}
}
\begin{tikzpicture}
  \def\vshift{0}
  \draw (-4.,1.7 + \vshift) node {\underline{Step 1}};
  \draw[<->] (-4.25,0 + \vshift) -- (4.25, 0 + \vshift) node [right]{};
  \draw[snake] (-3.5, 0 + \vshift) -- (.1,0 + \vshift) node [right]{};
  \draw (-3.5, 0.2cm + \vshift) -- (-3.5, -0.2cm + \vshift) (-3.5 cm,-3ex + \vshift) node
        [anchor=base,fill=white,inner sep=1pt]  {$-a$};
  \draw (0, 0.2cm + \vshift) -- (0,-0.2cm + \vshift) (0 cm,-3ex + \vshift) node
        [anchor=base,fill=white,inner sep=1pt]  {$0$};
  \draw (2, 0.2cm + \vshift) -- (2,-0.2cm + \vshift) (2 cm,-3ex + \vshift) node
        [anchor=base,fill=white,inner sep=1pt]  {$M_a$};
  \begin{scope}[shift={(-3.5, -0.2cm + \vshift)}]
    \lineann[1.]{0}{5.5}{$a + M_a$}
  \end{scope}
  \node [label={[align=left]If $a + M_a < L \implies$ switchback \\ If $a + M_a \geq L \implies$ no switchback}] at (0,-2. cm + \vshift){};

  \def\vshift{-4cm}
  \draw (-4.,1.7cm + \vshift) node {\underline{Step 2}};
  \draw[<->] (-4.25,0 + \vshift) -- (4.25, 0 + \vshift) node [right]{$(a'=M_a + a)$};
  \draw[snake] (-3.5, 0 + \vshift) -- (.1,0 + \vshift) node [right]{};
  \draw (-3.5, 0.2cm + \vshift) -- (-3.5, -0.2cm + \vshift) (-3.5 cm,-3ex + \vshift) node
        [anchor=base,fill=white,inner sep=1pt]  {$-a'$};
  \draw (0, 0.2cm + \vshift) -- (0,-0.2cm + \vshift) (0 cm,-3ex + \vshift) node
        [anchor=base,fill=white,inner sep=1pt]  {0};
  \draw (2, 0.2cm + \vshift) -- (2,-0.2cm + \vshift) (2 cm,-3ex + \vshift) node
        [anchor=base,fill=white,inner sep=1pt]  {$M_{a'}$};
  \begin{scope}[shift={(-3.5, -0.2cm + \vshift)}]
    \lineann[1.]{0}{5.5}{$a' + M_{a'}$}
  \end{scope}
  \node [label={[align=left]If $a' + M_{a'} < L \implies$ switchback \\ If $a' + M_{a'} \geq L \implies$ no additional switchback}] at (0,-2. cm + \vshift){};

  \def\vshift{-8cm}
  \draw (-4.,1.7cm + \vshift) node {\underline{Step 3}};
  \draw[<->] (-4.25,0 + \vshift) -- (4.25, 0 + \vshift) node [right]{$(a''=M_{a'} + a')$};
  \draw[snake] (-3.5, 0 + \vshift) -- (.1,0 + \vshift) node [right]{};
  \draw (-3.5, 0.2cm + \vshift) -- (-3.5, -0.2cm + \vshift) (-3.5 cm,-3ex + \vshift) node
        [anchor=base,fill=white,inner sep=1pt]  {$-a''$};
  \draw (0, 0.2cm + \vshift) -- (0,-0.2cm + \vshift) (0 cm,-3ex + \vshift) node
        [anchor=base,fill=white,inner sep=1pt]  {0};
  \draw (2, 0.2cm + \vshift) -- (2,-0.2cm + \vshift) (2 cm,-3ex + \vshift) node
        [anchor=base,fill=white,inner sep=1pt]  {$M_{a''}$};
  \begin{scope}[shift={(-3.5, -0.2cm + \vshift)}]
    \lineann[1.]{0}{5.5}{$a'' + M_{a''}$}
  \end{scope}
  \node [label={[align=left]If $a'' + M_{a''} < L \implies$ switchback \\ If $a'' + M_{a''} \geq L \implies$ no additional switchback}] at (0,-2. cm + \vshift){};
\end{tikzpicture}

%% file: BMCircle2.bbl
\begin{thebibliography}{00}

\bibitem[Borodin and Salminen(2002)]{Borodin}
Borodin, S. and Salminen, P.,
\emph{Handbook of Brownian Motion- Facts and Formulae}
Birkh{\"a}user Verlag (2002).

\bibitem[Feller (1951)]{Feller}
Feller, W.,
\emph{The Asymptotic Distribution of the Range of Sums of Independent Random Variables},
The Annals of Mathematical Statistics. \textbf{22} (1951), 427--432.

\bibitem[Imhof (1986)]{Imhof}
Imhof, J.P.,
\emph{On The Time Spent Above a Level by Brownian Motion with Negative Drift},
Advances in Applied Probability. \textbf{18} (1986), 1017--1018.


\bibitem[Oberhettinger and Badii (1973)]{Oberhettinger}
Oberhettinger, F. and Badii, L.,
\emph{Tables of Laplace Transforms}
Springer-Verlag (1973).

\bibitem[Vallois (1993)]{Vallois}
Vallois, P.,
\emph{Diffusion arr{\^e}t{\'e}e au premier instant o{\^u} l'amplitude atteint un niveau donn{\'e}},
Stochastics and Stochastic Reports. \textbf{43}, (1993) 93-115.







\end{thebibliography}
